\documentclass[10pt,conference]{IEEEtran}

\usepackage{graphics} 
\usepackage{epsfig} 
\usepackage{times} 
\usepackage{amsmath} 
\usepackage{amssymb}  
\usepackage{setspace}
\usepackage{color}
\usepackage{graphicx}

\newcommand{\mbf}[1]{\text{\boldmath{$#1$}}}
\newcommand{\by}{{\mbf{y}}}
\newcommand{\bY}{{\mbf{Y}}}

\newcommand{\mc}[1]{\mathcal{#1}} 
\newcommand{\cP}{{\mc{P}}}
\newcommand{\cY}{{\mc{Y}}}

\newcounter{thm}

\newcounter{rmk}
\newcounter{prp}

\newtheorem{remark}[rmk]{Remark}
\newtheorem{theorem}[thm]{Theorem}

\newtheorem{proposition}[prp]{Proposition}

\definecolor{dred}{rgb}{0,0,0}
\def\tdred{\textcolor{dred}}

\begin{document}

\title{Universal Sequential Outlier Hypothesis Testing}

\author{
\authorblockN{Yun Li\IEEEauthorrefmark{1}, Sirin Nitinawarat\IEEEauthorrefmark{2}, and Venugopal V. Veeravalli\IEEEauthorrefmark{1}}
\authorblockA{\IEEEauthorrefmark{1}
Department of Electrical and Computer Engineering\\
and\\
Coordinated Science Laboratory\\
University of Illinois at Urbana-Champaign\\ Urbana, IL 61801-2307, USA \\
Emails: \{yunli2, vvv\}@illinois.edu
}
\authorblockA{\IEEEauthorrefmark{2}
Qualcomm Technologies, Inc.\\
5775 Morehouse Drive\\
San Diego, CA 92121\\
Email: sirin.nitinawarat@gmail.com
}
}

\maketitle

\begin{abstract}
Universal outlier hypothesis testing is studied in a sequential setting.  Multiple observation sequences are collected, a small subset of which are outliers.  A sequence is considered an outlier if the observations in that sequence are generated by an ``outlier'' distribution, distinct from a common ``typical'' distribution governing the majority of the sequences.  Apart from being distinct, the outlier and typical distributions can be arbitrarily close.  The goal is to design a universal test to best discern all the outlier sequences.  A universal test with the flavor of the repeated significance test is proposed and its asymptotic performance is characterized under various universal settings.  The proposed test is shown to be universally consistent.  
\tdred{For the model with} identical outliers, the test is shown to be asymptotically optimal universally when the number of outliers is the largest possible and with the typical distribution being known, and its asymptotic performance otherwise is also characterized.  \tdred{An extension} of the findings to the model with multiple distinct outliers is \tdred{also discussed.}  In all cases, it is shown that the asymptotic performance guarantees for the proposed test when neither the outlier nor typical distribution is known converge to those when the typical distribution is known.
\end{abstract}

\section{Introduction} \label{sec-intro}

We consider the following inference problem of outlier hypothesis testing in a sequential setting.  Among a fixed number of independent and memoryless observation sequences, it is assumed that a small subset (possibly empty) of these sequences are outliers.  Specifically, most of the sequences are assumed to be distributed according to a ``typical'' distribution, while an outlier sequence is distributed according to an ``outlier distribution,'' distinct from the typical distribution.  We shall be interested in a {\em non-parametric} setting, in which the outlier and typical distributions are not fully known and can be arbitrarily close.  The goal is to design a universal test to identify all the outlier sequences using the fewest observations.

In \cite{li-niti-veer-ieeetit-2014}, we studied universal outlier hypothesis testing in a fixed sample size setting.  The main finding therein was that the generalized likelihood (GL) test is far more efficient for universal outlier hypothesis testing than for the other inference problems, such as homogeneity testing and classification \cite{pear-1911, ziv-1988, gutm-1989}.  In particular, the GL test was shown to be {\em universally exponentially consistent} for outlier hypothesis testing, whereas it is impossible to achieve universally exponential consistency for homogeneity testing or classification without training data \cite{ziv-1988, gutm-1989}.  We also showed that the GL test is {\em asymptotically optimal} in the limit of large number of sequences.  \tdred{Our previous paper \cite{li-niti-veer-isit-2014} generalized the scope of these previous findings to the sequential setting, but only covered
the setting with at most one outlier.  In this paper, we shall focus on the extension with multiple outliers.}  

Sequential hypothesis testing has a rich history going back to the seminal work of Wald \cite{wald-amstat-1945}.  A majority of the results on sequential hypothesis testing have been for the case where the conditional distributions of observations under the hypotheses are completely known (see, e.g., \cite{wald-amstat-1945, wald-wolf-amstat-1948, baum-veer-ieeetit-1994, veer-baum-ieeetit-1995, drag-tart-veer-ieeetit-1999, drag-tart-veer-ieeetit-2000}).  For the case where the distribution of observations is not completely specified, there have been a number of results for composite hypothesis testing with parametric families of distributions \cite{zack-1971, lai-astat-1988}.  As elucidated by Wald \cite{wald-amstat-1945}, there are two general approaches for constructing sequential tests for such parametric settings, one based on a weighted (or mixture) likelihood function for each hypothesis (see, e.g., \cite{zack-1971}), and the other based on a maximum (generalized) likelihood function for each hypothesis (see, e.g., \cite{lai-astat-1988}).  There have also been a limited number of papers on non-parametric approaches to sequential hypothesis testing where the functional form of the distribution is unknown, but it is known, for example, that the conditional distribution under the various hypotheses are rigid translations of each other (see, e.g., \cite{most-amstat-1948}).  Sequential outlier hypothesis testing is closely related to the so called {\em slippage problem} studied in the sequential setting (see, e.g., \cite{ferg-1967}).  In the slippage problem, $N$ populations are identically distributed except possibly for one.  The goal is to decide whether or not one of the populations has ``slipped'', if so, which one.  However, such prior work on the slippage problem concerned the situation where the typical and ``slipped'' distributions are tightly coupled, for example, when they are mean-shifted versions of each other. In universal sequential outlier hypothesis testing, we have no information regarding the outlier and typical distributions.  In particular, the typical and outlier distributions can be arbitrarily distributed and they can be arbitrarily close to each other.  In addition, we have no training data to learn these distributions before the test is performed. To the best of our knowledge, there has been no prior work on sequential outlier hypothesis testing in such a fully non-parametric setting that we study in this paper.  A key assumption that we make is that each instantaneous observation takes value in a finite and known set.  Under this assumption, we shall construct an efficient universal test that will be proven to be universally exponentially consistent, and to be sometimes optimal universally or in the limit of large number of sequences.  The proposed universal test has the flavor of the repeated significance test \cite{wood-nonlin-ren-th-book-1982, sieg-seq-anal-book-1985}, wherein the test stops when the generalized likelihood for the most likely hypothesis is larger by a time-dependent threshold than those for all the competing hypotheses, if that happens before a certain time.

\tdred{In Section \ref{sec-prelim}, we start with definitions of relevant distances between pairs of distributions, key to our results.}
\tdred{Sections \ref{sec-iden-outliers}, \ref{sec-distinct-outliers} concern the models with identical and distinctly distributed outliers, respectively.}  
\tdred{Performance of our proposed tests is evaluated on real data relevant to spam detection applications in Section \ref{sec-num-results}.}  
\tdred{Due to space limitations, proofs of our results are omitted.}

\section{PRELIMINARIES} 
\label{sec-prelim}

Throughout the paper, random variables (rvs) are denoted by capital letters, and their realizations are denoted by the corresponding lower-case letters.  All rvs are assumed to take values in {\em finite} sets, and all logarithms are the natural one.

Our results will be stated in terms of certain distance metrics between a pair of distributions $p, q$ on $\mathcal{Y}:$  the {\em Bhattacharyya distance} and the {\em relative entropy,} denoted by $B \left( p, q \right)$ and $D \left( p \| q \right),$ respectively, and defined as (see, e.g., \cite{cove-thom-eit-book-2006})
\begin{align}
B(p, q) \ \triangleq \ -\log \left ( \sum_{y \in \mathcal{Y}} p(y)^{\frac{1}{2}}q(y)^{\frac{1}{2}} \right ),  \nonumber
\end{align}
and
\begin{align}
D(p \| q) \ \triangleq \ \sum_{y \in \cY} p(y) \log \frac{p(y)}{q(y)}, \nonumber
\end{align}
respectively.

\section{MODEL WITH \tdred{IDENTICAL} OUTLIERS}
\label{sec-iden-outliers}

\tdred{Consider $M \geq 3$ independent sequences, each of which consists of independent and identically distributed (i.i.d.) observations.  Denote the $k$-th observation of the $i$-th sequence by $Y_k^{(i)} \in \mathcal{Y}$.  We} assume that there are {\em up to} $K > 2$ outliers among the $M$ sequences with $K < \frac{M}{2} $, and \tdred{that each of the outliers} are identically distributed (i.i.d.) according to \tdred{$\mu \in \cP \left( \cY \right)$, whereas all the other sequences are distributed according to the typical distribution $\pi \in \cP \left( \cY \right)$.  Under the hypothesis with no outlier, namely, the {\em null} hypothesis, all sequences are commonly distributed according to the typical distribution. {\em Nothing is known about $\mu$ and $\pi$ except that $\mu \neq \pi,$ and that each of them has a full support.}  The assumption of $\mu, \pi$ having full supports rules out trivial cases where it is easier to identify the outlier sequences.}
    
It was shown in \cite{li-niti-veer-ieeetit-2014} that in the fixed sample size setting, this assumption of the outliers being identically distributed is essential for the existence of a test that is universally exponentially consistent (under all the non-null hypotheses) when the number of outliers is not completely specified (anything from 0 to $K$).  In the next Section \ref{sec-distinct-outliers}, we shall look at the extension with possibly distinctly distributed outliers but with their total number being known.

\tdred{When there are some outliers, with the set of all outliers denoted by $S,\  0 < \vert S \vert < \frac{M}{2}$,}
the joint distribution of \tdred{the first $n$ observations} is given by
\begin{align}
p_S \left( \by^{n} \right)
&\ =\ p_S \left( \by_{1}, \ldots, \by_{n} \right)	\nonumber \\
&\ =\ \prod_{k=1}^n 
\left \{ 
\prod_{i \in S} \mu \left( y_k^{(i)} \right)
\prod_{j \notin S} \pi \left( y_k^{(i)} \right)
\right \}.
\label{eqn-likelihood-commom-mu}
\end{align}
\tdred{Under the null hypothesis with no outlier, the joint distribution of the observations is given as
\begin{align}
p_{0} \left ( \by^{n} \right ) 
\ =\ p_{\emptyset} \left( \by^n \right)
\ = \ \prod_{k=1}^{n} \prod_{i=1}^{M} \pi \left ( y_{k}^{(i)} \right ).
\nonumber
\end{align}
}

\tdred{A sequential test for the outlier consists of a stopping rule and a final decision rule.  The stopping rule defines a random (Markov) time, denoted by $N$, which is the number of observations taken until a final decision is made.  At the stopping time $N = n,$ a decision is made based on a decision rule $\delta: \cY^{Mn} \rightarrow \mathcal{S},$ where $\mathcal{S}$ denotes the set of all subsets of $\left \{1, \ldots, M \right \}$ of size at most $K$.  The overall goal of sequential testing is to achieve a certain level of accuracy for the final decision using the fewest number of observations on average.}

\tdred{
We consider the sequential outlier hypothesis testing problem in two settings: the setting where only $\pi$ is known, and the completely universal setting where neither $\mu$ nor $\pi$ is known.  Consequently, a universal test is not allowed to be a function of $\mu$, and of $\mu$ or $\pi$, in the respective settings.
}

\tdred{
The accuracy of a sequential test is gauged using the maximal error probability $P_{\max}$, which is a function of both the test and $\left( \mu, \pi \right)$ and is defined as
\begin{align}
P_{\max} &\ \triangleq \ \max_{S \in \mathcal{S}} \ 
\mathbb{P}_{S} \left \{ \delta \left( \bY^{N} \right) \neq S \right \},  
\label{eqn-def-maxerr}
\end{align}
where $\mathbb{P}_S,\ S \in \mathcal{S},$ denotes the probability distribution for the hypothesis with $S$ as the subset of all outliers.  We say a sequence of tests is {\em universally consistent} if the maximal error probability converges to zero for any $\mu, \pi, \mu \neq \pi$.   Further, we say it is {\em universally exponentially consistent} if the exponent for the maximal error probability with respect to the expected stopping time under each hypothesis
is strictly positive, i.e., there exists $\alpha_S > 0$ such that 
\begin{align}
\mathbb{E}_S \left [ N \right ]
~\leq~
\frac{- \log P_{\max}}{\alpha_S} \left( 1 + o (1) \right) \label{eqn-def-uec}
\end{align}
for any $\mu, \pi, \mu \neq \pi$ as $P_{\max} \rightarrow 0$.
}

\tdred{
We first consider the setting where both the typical and outlier distributions are known.  In this non-universal setting, the Multihypothesis Sequential Probability Ratio Test (MSPRT) 
was shown to be asymptotically optimal in the regime with vanishing error probability \cite{drag-tart-veer-ieeetit-1999}.  For a given threshold $T > 0$ and with $\hat{S}  \left( \by^{n} \right)~\triangleq~\mathop{\mathrm{argmax}}\limits_{S \in \mathcal{S}} p_S \left ( \by^{n} \right ),$ denoting the instantaneous maximum likelihood (ML) estimate of the hypothesis at time $n$, the stopping time $N^*$ and final decision $\delta^*$ of the MSPRT are defined as follows.
\begin{align}
N^* & \  = \ 
\mathop{\mbox{argmin}}_{n \geq 1} \ 
\left [
\frac{p_{\hat{S}} \, \left( \bY^{n} \right)}
{\max\limits_{S \neq \hat{S}} \ p_{S} \left( \bY^{n} \right)} 
\ > \ T
\right ], 
\label{eqn-MSPRTstoprule-mult-outlier} \\
\delta^* & \ =\  \hat{S} \left( \bY^{N^*} \right).  
\label{eqn-MSPRTdecision-mult-outlier}
\end{align}
}


\begin{proposition}
As the threshold $T$ in (\ref{eqn-MSPRTstoprule-mult-outlier}) approaches infinity, the MSPRT in (\ref{eqn-MSPRTstoprule-mult-outlier}), (\ref{eqn-MSPRTdecision-mult-outlier}) satisfies
$P_{\max} \ =\ O \left( \frac{1}{T} \right).$ 
In addition, it \tdred{yields that
\begin{align}
\mathbb{E}_{S} \left [ N^* \right ] =
\left \{
\begin{array}{cc}
\frac{-\log{P_{\max}}}{D \left( \mu \| \pi \right) } (1+ o(1)), 
& \vert S \vert = K;\\
\frac{-\log{P_{\max}}}{\min \left ( D \left( \mu \| \pi \right), D \left( \pi \| \mu \right)  \right ) } (1+ o(1)), 
& 1 \leq \vert S \vert < K;\\
\frac{-\log{P_{\max}}}{D \left( \pi \| \mu \right) } (1+ o(1)),
& S = \emptyset.
\end{array}
\right. 
\nonumber
\end{align}
}

Furthermore, the MSPRT is asymptotically optimal. In particular, for any sequence of tests $\left( N, \delta \right)$ with vanishing maximal error probability, it \tdred{holds (simultaneously) that
\begin{align}
\mathbb{E}_{S}[N] \geq
\left \{
\begin{array}{cc}
\frac{-\log P_{\max}}{D \left( \mu \| \pi \right) } (1+ o(1))
& \vert S \vert = K; \\
\frac{-\log P_{\max}}{\min \left \{ D \left( \mu \| \pi \right), D \left( \pi \| \mu \right)  \right \} } (1+ o(1))
& 1 \leq \vert S \vert < K; \\
\frac{-\log P_{\max}}{D \left( \pi \| \mu \right) } (1+ o(1))
& S = \emptyset.
\end{array}
\right.
\nonumber
\end{align}
}
\end{proposition} 

\tdred{Now we consider the universal settings when the outlier distribution is unknown, and when neither the outlier nor typical distribution is known.  In the fixed sample size setting, it was shown in \cite{li-niti-veer-ieeetit-2014} that a universally exponentially consistent test cannot exist.  Correspondingly, we proposed a test therein that fulfilled  a lesser objective of attaining universally exponential consistency under all the non-null hypotheses, while satisfying {\em only} universal consistency under the null hypothesis.  We now describe a universal sequential test fulfilling a similar objective.}

\subsubsection{Proposed Universal Test}

\tdred{For each $i = 1, \ldots, M$, denote the empirical distribution of $\by^{(i)}$ by $\gamma_i.$}  When only $\pi$ is known, we compute the generalized likelihood of $\by^{n}$ under each non-null hypothesis corresponding to a non-empty subset $S \subset \left \{ 1, \ldots, M \right \}$ by replacing the unknown $\mu$ in (\ref{eqn-likelihood-commom-mu}) with its ML estimate $\hat{\mu}_{S} \triangleq \frac{\sum_{i \in S} \gamma_{i}}{\vert S \vert}$, as 
\begin{align}
\hat{p}^{\text{typ}}_{S} \left ( \by^{n} \right ) 
& \ = \ 
\prod_{k=1}^{n} 
\left \{ 
\prod_{i \in S}  \hat{\mu}_{S} \left ( y_{k}^{(i)} \right ) \, 
\prod_{j \notin S} \pi \left ( y_{k}^{(j)} \right ) \right \}. 
\label{eqn-GL-mult-outlier-piknown}
\end{align}
Similarly, when neither $\pi$ nor $\mu$ is known, we compute the generalized likelihood of $\by^{n}$ under each non-null hypothesis corresponding to a non-empty $S \in \mathcal{S}$ by replacing the unknown $\mu$ and $\pi$ in (\ref{eqn-likelihood-commom-mu}) with their ML estimates $\hat{\mu}_{S} \triangleq \frac{\sum_{i \in S} \gamma_{i}}{\vert S \vert}$, and $\hat{\pi}_{S} \triangleq 
\frac{\sum_{j \notin S} \gamma_{j}}{M-\vert S \vert}$, respectively, as 
\begin{align}
\hat{p}^{\text{univ}}_{S} \left ( \by^{n} \right ) 
& \ = \ 
\prod_{k=1}^{n} 
\left \{ 
\prod_{i \in S} \hat{\mu}_{S} \left ( y_{k}^{(i)} \right ) \, 
\prod_{j \notin S} \hat{\pi}_{S} \left ( y_{k}^{(j)} \right ) \right \}. 
\label{eqn-GL-mult-outlier-piunknown}
\end{align}

When only $\pi$ is known and with $\hat{S} \left( \bY^{N} \right) ~\triangleq$\\
$\mathop{\mathop{\mathrm{argmax}}\limits_{S \in \mathcal{S}}}\limits_{S \neq \emptyset} 
\hat{p}^{\text{typ}}_{S} \left ( \by^{n} \right )   
= \mathop{\mathop{\mathrm{argmin}}\limits_{S \in \mathcal{S}}}\limits_{S \neq \emptyset}
 \left [
\sum\limits_{i \in S} D \Big (\gamma_{i} \| 
{\textstyle \frac{\sum_{k \in S} \gamma_{k}}{\vert S \vert} \Big )}
+ \sum\limits_{j \notin S} D(\gamma_{j} \| \pi) 
\right ]$ denoting the instantaneous estimate of the non-null hypothesis (using the generalized likelihood) at time $n$, our proposed universal test can be described by the following stopping and final decision rules
\begin{align}
N^* &\ = \ \min \left ( \tilde{N}, \left \lfloor \tdred{f(T)}  \right \rfloor  \right ), 
\label{eqn-stoprule-unknown-mult-outlier-known} \\
\delta^* &\ =\ 
\left \{
\begin{array}{cc}
\hat{S} \left( \bY^{N^*} \right)		&\mbox{if}~\tilde{N} \leq \tdred{f(T)}; \\
0	&\mbox{if}~\tilde{N} > \tdred{f(T)},
\end{array}
\right.
\label{eqn-decision-unknown-mult-outlier-known-no-1}
\end{align}
where \tdred{$f(T)$ is any function growing at least as fast as $T \log{T}$, and}
\begin{align}
\tilde{N} &\triangleq \mathop{\mbox{argmin}}_{n \geq 1} 
\Bigg [
\mathop{\min\limits_{S' \neq \hat{S}}}\limits_{S' \neq \emptyset}
n 
\Big [ 
 \sum\limits_{i \in S'} 
D \left (\gamma_{i} \big \| \textstyle{\frac{\sum_{k \in S'} \gamma_{k} }{\vert S' \vert }} 
\right) + \sum\limits_{j \notin S'} 
D \left (\gamma_{j} \big \| \pi  \right)   
\nonumber \\
& \hspace{1.1in}
- \sum\limits_{i \in \hat{S}} 
D \left (\gamma_{i} \big \|  \textstyle{\frac{\sum_{k \in \hat{S}} \gamma_{k}}{\vert S \vert }}  
\right)
- \sum\limits_{j \notin \hat{S}} 
D \left (\gamma_{j} \big \| \pi  \right) 
\Big ] 
\nonumber \\
& \hspace{0.8in} > \log{T} + (M+1) \vert \mathcal{Y} \vert \log (n+1)
\Bigg ], 
\label{eqn-stoprule-unknown-mult-outlier-known-no-1}  
\end{align}

Similarly, when neither $\mu$ nor $\pi$ is known, the test can be written as
in (\ref{eqn-stoprule-unknown-mult-outlier-known}), (\ref{eqn-decision-unknown-mult-outlier-known-no-1}) but with $\hat{S} \left( \bY^{N} \right)
\triangleq \mathop{\mathrm{argmax}}\limits_{S \in \mathcal{S}, S \neq \emptyset} 
\hat{p}^{\text{univ}}_{S} \left ( \by^{n} \right )   
= \mathop{\mathrm{argmin}}\limits_{S \in \mathcal{S}, S \neq \emptyset} 
\Big [ 
\sum\limits_{i \in S} 
D \Big (\gamma_{i} \| {\textstyle  \frac{\sum_{k \in S} \gamma_{k}}{\vert S \vert}} \Big ) 
+
\sum\limits_{j \notin S} 
D \Big (\gamma_{j} \| {\textstyle \frac{\sum_{k \notin S} \gamma_{k} }{M- \vert S \vert }} \Big )
\Big ]$, and
\begin{align}
\tilde{N} &\triangleq \mathop{\mbox{argmin}}_{n \geq 1}
\Bigg [ 
\mathop{\min\limits_{S' \neq \hat{S}}}\limits_{S' \neq \emptyset}\ 
n 
\Big [
 \sum\limits_{i \in S'} 
D \left (\gamma_{i} \big \| \textstyle{\frac{\sum_{k \in S'} \gamma_{k} }{\vert S' \vert }} 
\right) 
\nonumber \\
& \hspace{1.3in} + 
\sum\limits_{j \notin S'} 
D \left (\gamma_{j} \big \| \textstyle{\frac{\sum_{k \notin S'} \gamma_{k} }{M-\vert S' \vert }} 
\right)  
\nonumber \\
& \hspace{1.3in} 
- \sum\limits_{i \in \hat{S}} 
D \left (\gamma_{i} \big \|  \textstyle{\frac{\sum_{k \in \hat{S}} \gamma_{k}}{\vert S \vert }}  
\right)
\nonumber \\
& \hspace{1.3in} 
- \sum\limits_{j \notin \hat{S}} 
D \left (\gamma_{j} \big \|  \textstyle{\frac{\sum_{k \notin \hat{S}} \gamma_{k}}{M- \vert S \vert}}  
\right) \Big ]
\nonumber \\
& \hspace{0.8in} > \log{T} + (M+1) \vert \mathcal{Y} \vert \log (n+1)
\Bigg ].
\label{eqn-stoprule-unknown-mult-outlier-known-no-2}
\end{align}

\subsubsection{Performance of Proposed Test}

\setcounter{thm}{1}
\begin{theorem} \label{thm-5}
When only $\pi$ is known, the test in (\ref{eqn-stoprule-unknown-mult-outlier-known}), (\ref{eqn-decision-unknown-mult-outlier-known-no-1}), (\ref{eqn-stoprule-unknown-mult-outlier-known-no-1}) is universally consistent, and yields for every $T$ that
$P_{\max} \ \leq \ \frac{C}{T},$  
where $C$ is a constant that depends only on $M, K, \mu, \pi$, but not on $T$.
In addition, it satisfies for each non-null hypothesis $S \in \mathcal{S}, S \neq \emptyset,$ that as $T \rightarrow \infty$,
\begin{align}
\mathbb{E}_{S}[N^*] &\ \tdred{\leq} \ \frac{\log{T}}{\alpha_{S}}(1+o(1))   
\nonumber \\
&\ \leq\  
\left \{
\begin{array}{cc}
\frac{-\log{P_{\max}}}{D(\mu \| \pi)}(1+o(1)), & \vert S \vert = K; \\
\frac{-\log{P_{\max}}}
{\min \left( D(\mu \| \pi), \eta_{S}(\mu \| \pi)\right)}(1+o(1)), & 1 \leq \vert S \vert < K,
\end{array}
\right.	
\label{eqn-lowerbd-alphaSb}
\end{align}
where 
\begin{align}
\alpha_{S} &\triangleq 
\mathop{\min\limits_{S' \neq S}}\limits_{S' \neq \emptyset}
\bigg [ \  
 \vert S \cap S' \vert 
 D \left ( \mu \Big \| \ 
\textstyle{\frac{ \vert S \cap S' \vert \mu + \vert S' \backslash S \vert \pi}
{\vert S' \vert }} 
\right ) 
\nonumber \\
& \hspace{0.55in}
+ 
\big \vert S \backslash S' \big  \vert D (\mu  \| \pi  ) 
 \nonumber \\
& \hspace{0.55in}
+  
\big \vert S' \backslash S \big \vert D \left ( \pi \Big  \| \ 
{\textstyle  
\frac{ \vert S \cap S' \vert \mu + \vert S' \backslash S \vert \pi}{\vert S' \vert } }
\right) 
\  \bigg ]
\ > \ 0. 
\nonumber
\end{align}
and
\begin{align}
\eta_{S}(\mu \| \pi)  \ \triangleq \ \min_{p \in \mathcal{P}(\mathcal{Y})} \ \vert S \vert D(\mu \| p) + D(\pi \| p). \label{eqn-eta}
\end{align}
\end{theorem}

\begin{theorem} \label{thm-6}
When neither $\mu$ nor $\pi$ is known, the universal test in (\ref{eqn-stoprule-unknown-mult-outlier-known}), (\ref{eqn-decision-unknown-mult-outlier-known-no-1}), (\ref{eqn-stoprule-unknown-mult-outlier-known-no-2}) is universally consistent, and yields for every $T$ that
$P_{\max} \ \leq \ \frac{C}{T},$  
where $C$ is a constant that depends on $M, K, \mu, \pi$, but not on $T$.
In addition, for each non-null hypothesis $S \in \mathcal{S}, S \neq \emptyset,$ as $T \rightarrow \infty$,
\begin{align}
\mathbb{E}_{S} \left [ N^* \right] &\ \tdred{\leq} \ \frac{\log T}{\overline{\alpha}_{S}}(1+o(1))   
\label{eqn-thm7-claim2} \\
&\ \leq\ \left \{
\begin{array}{cc}
\frac{-\log{P_{\max}}}{ \overline{\eta}(\mu \| \pi)}(1+o(1)), & \vert S \vert = K;\\
\frac{-\log{P_{\max}}}
{\min \left( \overline{\eta}(\mu \| \pi), \eta_{S}(\mu \| \pi)\right)}(1+o(1)), 
& 1 \leq \vert S \vert < K,
\end{array}
\right.
\label{eqn-lowerbd-alphaS-1b} 
\end{align} 
where
\begin{align}
\hspace{-0.1in}
\overline{\alpha}_{S} &\triangleq
\mathop{\min\limits_{S' \neq S}}\limits_{S' \neq \emptyset}
\bigg [ 
\big \vert S \cap S' \big \vert 
D \Big ( \mu \Big \| \ 
{\textstyle  \frac{ \vert S \cap S' \vert \mu + \vert S' \backslash S \vert \pi}{\vert S' \vert } }
\Big )
\nonumber \\
& \hspace{0.5in} 
+ 
\big \vert S \backslash S' \big \vert 
D 
\Big (\mu \Big \| \ 
{\textstyle  
\frac{ \vert S \backslash S' \vert \mu + \vert S^c \cap S'^c \vert \pi}{M- \vert S' \vert } 
}
\Big )  
\nonumber \\
&\hspace{0.5in} 
+
\big \vert S' \backslash S \big \vert 
D \Big (\pi \Big  \| \ 
{\textstyle  
\frac{ \vert  S \cap S' \vert \mu + \vert S' \backslash S \vert \pi}{\vert S' \vert } 
}
\Big )
\nonumber \\
&\hspace{0.5in} 
+ 
\big \vert  S^c \cap S'^{c} \big \vert 
D \Big (\pi \Big  \| \ 
{\textstyle  
\frac{ \vert S \backslash S' \vert \mu + \vert S^c \cap S'^c  \vert \pi}{M- \vert S' \vert } 
}
\Big )
\bigg ] \ > \ 0, 
\nonumber
\end{align}
and
\begin{align}
\overline{\eta}_{S} \left( \mu \| \pi \right) 
\triangleq \min_{p \in \mathcal{P}(\mathcal{Y})} D(\mu \| p) + (M - K - \vert S \vert) D(\pi \| p),
\label{eqn-eta-bar}
\end{align}
and $\eta_{S} \left( \mu \| \pi \right)$ is as in (\ref{eqn-eta}).
\end{theorem}
\begin{remark}
It follows from (\ref{eqn-eta-bar}) that as $M \rightarrow \infty,$
\begin{align}
\overline{\eta}_{S} \left( \mu, \pi \right)
\ \rightarrow\  D(\mu \| \pi), \label{eqn-prop8-claim2}
\end{align}
i.e., the asymptotic performance guarantee for the test in (\ref{eqn-stoprule-unknown-mult-outlier-known}), (\ref{eqn-decision-unknown-mult-outlier-known-no-1}), (\ref{eqn-stoprule-unknown-mult-outlier-known-no-2}) when neither $\mu$ nor $\pi$ (cf. (\ref{eqn-lowerbd-alphaS-1b})) are known converges to that for the test in (\ref{eqn-stoprule-unknown-mult-outlier-known}), (\ref{eqn-decision-unknown-mult-outlier-known-no-1}), (\ref{eqn-stoprule-unknown-mult-outlier-known-no-1}) when  $\pi$ is known (cf. (\ref{eqn-lowerbd-alphaSb})) as $M \rightarrow \infty.$
\end{remark}

\section{\tdred{MODEL WITH DISTINCT OUTLIERS}}
\label{sec-distinct-outliers}

It was shown in \cite{li-niti-veer-ieeetit-2014} that when the outliers can be arbitrarily distinctly distributed, the assumption of the number of outliers being known is essential for the existence of a universally exponentially consistent test.  We now describe this extension with distinctly distributed outliers but with their number being known in the sequential setting.  

In particular, for an $S \subset \left \{1, \ldots, M \right \},\ \vert S \vert = K,$ denoting the set of $K$ outliers, the joint distribution of all observations under the hypothesis with the outlier subset being $S$ is 
\begin{align}
p_S \left( \by^{n} \right)
&\ =\ p_S \left( \by_{1}, \ldots, \by_{n} \right)	\nonumber \\
&\ =\ \prod_{k=1}^n 
\left \{ 
\prod_{i \in S} \mu_i \left( y_k^{(i)} \right)
\prod_{j \notin S} \pi \left( y_k^{(i)} \right)
\right \},
\label{eqn-likelihood-distinct-mu}
\end{align}
where each $i$-th outlier, $i \in S,$ is distributed as $\mu_i,$ which can be arbitrarily distinct from one another as long as each $\mu_i \neq \pi.$  At the stopping time $N = n,$ the test for the outliers is done based on a rule $\delta: \mathcal{Y}^{Mn} \rightarrow \mathcal{S}_K,$ where $\mathcal{S}_K$ will now denote the set of all subsets of $\left \{ 1, \ldots, M \right \}$ of size {\em exactly} $K$.  Notice that unlike in the previous sections, the current model does not include the null hypothesis with no outlier.  The maximal error probability is defined as previously in (\ref{eqn-def-maxerr}) but with the maximum being over $\mathcal{S}_K$ instead.

\subsection{Proposed Universal Test}

When only $\pi$ is known, we can compute the corresponding generalized likelihood of $\by^{n}$ under each hypothesis $S \in \mathcal{S}_K$ by replacing the unknown $\mu_i,\ i \in S,$ in (\ref{eqn-likelihood-distinct-mu}) with its ML estimate $\hat{\mu}_{S}^i \triangleq \gamma_i.$  In particular, with $\hat{S} \left( \bY^{n} \right)
=  \mathop{\mathrm{argmin}}\limits_{S \in \mathcal{S}_K}
\sum\limits_{j \notin S} D(\gamma_{j} \| \pi)$ denoting the instantaneous estimate of the hypothesis (using the generalized likelihood) at time $n$, our proposed universal test can be described by the following stopping and final decision rules: 
\begin{align}
N^* &\ =\ \mathop{\mbox{argmin}}_{n \geq 1} 
\Bigg [
\mathop{\min_{S' \neq \hat{S}}}_{S' \in \mathcal{S}_K} \ n 
\Big [ 
\sum\limits_{j \notin S'} 
D \left (\gamma_{j} \big \| \pi  \right) 
- \sum\limits_{j \notin \hat{S}} 
D \left (\gamma_{j} \big \| \pi  \right) 
\Big ]	
\nonumber \\
& \hspace{0.75in}
> \log{T} + (M+1) \vert \mathcal{Y} \vert \log (n+1)
\Bigg ];
\label{eqn-stoprule-known-mult-outlier-knownmu}  \\
\delta^* &\ =\  \hat{S} \left( \bY^{N^*} \right).
\label{eqn-decision-known-mult-outlier-knownmu}
\end{align}

Similarly, when neither $\mu$ nor $\pi$ is known, the test can be written as
\begin{align}
N^* &= 
\mathop{\mbox{argmin}}_{n \geq 1} \, 
\Bigg [ 
\mathop{\min_{S' \neq \hat{S}}}_{S' \in \mathcal{S}_K} \ 
n 
\Big [
\sum\limits_{j \notin S'} 
D \left (\gamma_{j} \big \| \textstyle{\frac{\sum_{k \notin S'} \gamma_{k} }{M-\vert S' \vert }} 
\right) 
\nonumber \\
& \hspace{1.3in}
- \sum\limits_{j \notin \hat{S}} 
D \left (\gamma_{j} \big \|  \textstyle{\frac{\sum_{k \notin \hat{S}} \gamma_{k}}{M- \vert S \vert}}  
\right) \Big ]
\nonumber \\
& \hspace{0.8in}
> \log{T} + (M+1) \vert \mathcal{Y} \vert \log (n+1)
\Bigg ];
\label{eqn-stoprule-known-mult-outlier-unknownmu} \\
\delta^* &\ =\  \hat{S} \left( \bY^{N^*} \right),
\label{eqn-decision-known-mult-outlier-unknownmu}
\end{align}
but with 
$\hat{S} \left( \bY^{n} \right)
= \mathop{\mathrm{argmin}}\limits_{S \in \mathcal{S}} 
\sum\limits_{j \notin S} 
D \Big (\gamma_{j} \| {\textstyle \frac{\sum_{k \notin S} \gamma_{k} }{M- \vert S \vert }} \Big ).$
Note that since the null hypothesis is not present in this case, there is no need to truncate 
the stopping time by a predefined horizon as in (\ref{eqn-stoprule-unknown-mult-outlier-known}).

\subsection{Performance of the Proposed Tests}


\begin{theorem} 
\label{thm-7}
With the number of distinct outliers $K$ being known and when only $\pi$ is known, the test in (\ref{eqn-stoprule-known-mult-outlier-knownmu}), (\ref{eqn-decision-known-mult-outlier-knownmu}) is universally exponentially consistent, and yields for every $T$ that
$P_{\max} \ \leq \ \frac{C}{T},$ 
where $C$ is a constant that depends only on $M, K, \mu, \pi$, but not on $T$.
In addition, for each hypothesis $S \in \mathcal{S}_K$, as $T \rightarrow \infty$,
\begin{align}
\mathbb{E}_{S} \left[ N^* \right] &\ \leq \ 
\frac{-\log{P_{\max}}}{\min\limits_{i \in S} D \left( \mu_i \| \pi \right)}(1+o(1)).
\label{eqn-thm7-claim}
\end{align}
\end{theorem}

\begin{theorem} 
\label{thm-8}
With the number of distinct outliers $K$ being known, but neither $\mu$ nor $\pi$ being known, the test in (\ref{eqn-stoprule-known-mult-outlier-unknownmu}), (\ref{eqn-decision-known-mult-outlier-unknownmu}) is universally exponentially consistent, and yields for every $T$ that
$P_{\max} \ \leq \ \frac{C}{T},$
where $C$ is a constant that depends only on $M, K, \mu, \pi$, but not on $T$.
In addition, for each hypothesis $S \in \mathcal{S}_K$, as $T \rightarrow \infty$,
\begin{align}
\mathbb{E}_{S} \left[ N^* \right]  \leq 
\frac{-\log{P_{\max}}(1+o(1))}
{
\min\limits_{i \in S}
\min\limits_{p} 
\left( 
D \left( \mu_i \| p \right)
+ (M-2K) D \left( \pi \| p \right)
\right)
}.
\label{eqn-thm8-claim}
\end{align}
\end{theorem}
\begin{remark}
As $M \rightarrow \infty,$ the inner minimum in the denominator in (\ref{eqn-thm8-claim}) is attained at $p^* = \pi$ and, hence, the coefficient therein converges to $\min\limits_{i \in S} D \left( \mu_i \| p \right),$ which is the asymptotic performance of the universal test in 
(\ref{eqn-stoprule-known-mult-outlier-knownmu}), (\ref{eqn-decision-known-mult-outlier-knownmu}) when $\pi$ is known (cf. (\ref{eqn-thm7-claim})).
\end{remark}

\section{APPLICATION TO SPAM DETECTION}
\label{sec-num-results}

We evaluate the performance of the proposed universal tests on real data set relevant to spam detection applications.  The labeled data set (spam or non-spam) contains information from 4610 emails addressed to an employee at Hewlett-Packard and is publicly available \cite{hast-tibs-frie-2009}.  In particular, the data set consists of relative frequencies of a set of 48 words and 6 punctuation marks.  Out of 4601 emails, there are 1813 spams.

We design an experiment for the case with $M = 5$ sequences, and with at most two identical outliers for the total number hypotheses of 16.  Instead of looking at the frequencies of all words and punctuation marks available in the data set, we just pick out the frequency of the word ``HP,'' and quantize it into 5 levels (in the original data set, the frequencies take continuous values in the interval $[0, 100]$).  The $f(T)$ in (\ref{eqn-stoprule-unknown-mult-outlier-known}), (\ref{eqn-decision-unknown-mult-outlier-known-no-1}) is selected to be $T^5.$  The values of $\frac{- \log{P_{\max}}}{\mathbb{E}_S \left [ N^* \right]},$ for $\vert S \vert = 1$, and $\vert S \vert = 2,$ achievable by the universal test in  (\ref{eqn-stoprule-unknown-mult-outlier-known}), (\ref{eqn-decision-unknown-mult-outlier-known-no-1}), (\ref{eqn-stoprule-unknown-mult-outlier-known-no-2}) (without knowledge of either $\mu$ or $\pi$) is listed as a function of $T$ in Table I.  The asymptote of the expected stopping time (relative to the exponent for the error probability) under a hypothesis with $\vert S \vert = 2$ is lower than that under a hypothesis with $\vert S \vert = 1,$ which agrees with the results in (\ref{eqn-thm7-claim2}) and (\ref{eqn-lowerbd-alphaS-1b}).

\begin{table}[h]
	\caption[Table 1]{}
	\begin{tabular}{| c || c | c | c | c |}
	\hline
		 & $T = 3.98$ & $T = 4$ & $T = 4.05$ & $T = 4.1$ \\
	\hline 
	\hline
	${\textstyle \frac{- \log{P_{\max}}}{\mathbb{E}_{ \{ 1 \} } \left[ N^* \right]}}$
		& 0.0012 & 0.0017 & 0.0039 & 0.0069 \\
		& & & &\\
	\hline
	${\textstyle \frac{- \log{P_{\max}}}{\mathbb{E}_{ \{ 1, 2 \} } \left[ N^* \right]}}$
		& 0.0017 & 0.0025 & 0.0057 & 0.01 \\
		& & & &\\
	\hline
	\end{tabular}
\end{table}
 
\section*{ACKNOWLEDGMENT}

This research was partially supported by the Air Force Office of Scientific Research (AFOSR) under Grant FA9550-10-1-0458 through the University of Illinois at Urbana-Champaign, and by the National Science Foundation under Grant NSF CCF 11-11342.

\bibliographystyle{IEEEtran}
\bibliography{li-niti-veer-asilomar-2014}

\end{document}